\documentclass[10pt,reqno]{amsart}
\usepackage{amssymb}

\theoremstyle{plain}

\theoremstyle{remark}

\theoremstyle{definition}

\numberwithin{equation}{section}

\newcommand{\eqto}{\stackrel{\lower1.5pt\hbox{$\scriptstyle\sim\,$}}\to}

\def\Spec{{\rm Spec}}

\def\inv{{\rm inv}}

\def\A{{\mathbb A}}

\def\PP{{\mathbb P}}
\def\Q{{\mathbb Q}}

\def\Z{{\mathbb Z}}

\def\Br{{\rm Br}}

\begin{document}
\title[Brauer--Manin obstructions]
{Brauer--Manin obstructions to integral points}
\author{Andrew Kresch}
\address{
  Institut f\"ur Mathematik,
  Universit\"at Z\"urich,
  Winterthurerstrasse 190,
  CH-8057 Z\"urich, Switzerland
}
\email{andrew.kresch@math.unizh.ch}
\author{Yuri Tschinkel}
\address{
  Mathematisches Institut,
  Georg-August-Universit\"at G\"ottingen,
  Bunsenstrasse 3-5,
  D-37073 G\"ottingen, Germany
  and
  Courant Institute,
  251 Mercer Street,
  New York, NY 10012
}
\email{tschinkel@cims.nyu.edu}

\date{7 September 2007}
\subjclass[2000]{14G25 (primary); 14F22 (secondary).}
\thanks{The second author was supported by the NSF}

\begin{abstract}
We study Brauer--Manin obstructions to integral points on
open subsets of the projective plane.
\end{abstract}
\maketitle

\section{Introduction}
\label{sec:introduction}

Let $k$ be a number field and $X$ a smooth projective
geometrically irreducible variety
over $k$.
It is well known that
the existence of points of $X$ over all completions $k_v$ of $k$
does not imply the existence of a $k$-rational point on $X$,
in general.
This phenomenon is referred to as the failure of the
Hasse principle.
Examples of failure of the Hasse principle are known
for genus 1 curves, cubic surfaces, etc.
Even when the Hasse principle holds,
rational points need not be dense in the
set of adelic points of $X$.
This phenomenon, the failure of weak approximation,
also is known in many examples.
For instance, weak approximation always fails when
$X(k)\ne \emptyset$ and
$X\otimes \bar k$ (the variety obtained by base change to the
algebraic closure $\bar k$ of $k$)
has nontrivial algebraic fundamental group \cite{minchev}.

Often (e.g., conjecturally for rational surfaces)
the failure of the Hasse principle and weak approximation
is explained by the Brauer--Manin obstruction \cite{manin}, \cite{maninbook}.
By the exact sequence from class field theory
$$0\to \Br(k)\longrightarrow \bigoplus_v \Br(k_v)
\stackrel{\sum \inv_v}\longrightarrow \Q/\Z\to 0$$
(where $\sum \inv_v$ denotes the sum of local invariants)
we have the constraint
$$X(k)\subset X(\A_k)^{\Br}:=\{\,(x_v)\in X(\A_k)\,|\,
\sum \inv_v(\alpha|_{x_v})=0\,\forall\,\alpha\in\Br(X)\,\}$$
on the set $X(\A_k)$ of adelic points on $X$.
When $X(\A_k)^{\Br}\neq X(\A_k)$, then we say there is a
Brauer--Manin obstruction to the Hasse principle,
respectively to weak approximation, in case
$X(\A_k)^{\Br}=\emptyset$, respectively,
$X(\A_k)^{\Br}\neq\emptyset$.

For a thorough introduction to the subject, see
\cite{skorobogatovbook}.
For a survey, see \cite{peyre}.

The study of rational points on projective hypersurfaces is equivalent
to the study of integral solutions to
homogeneous Diophantine equations $f(x_0,\ldots,x_n)=0$.
Many interesting Diophantine problems involve
non-homogeneous equations.
Their solutions can be interpreted as integral points
on quasi-projective varieties.

Let $\mathfrak{o}_k$ be the ring of integers of $k$.
Let $\mathcal{X}$ be an integral model for $X$, i.e., a scheme,
projective and flat over $\Spec(\mathfrak{o}_k)$ having general fiber $X$.
Let $Z$ be a reduced closed subscheme of $X$, and set
$U=X\smallsetminus Z$.
Let $\mathcal{Z}$
be an integral model for $Z$
(which is uniquely determined by $Z$ and $\mathcal{X}$).
An \emph{integral point} of $U$ is a $k$-rational point on $X$ whose
unique extension to an $\mathfrak{o}_k$-point of $\mathcal{X}$ has
image disjoint from $\mathcal{Z}$.
One can also speak of $S$-integral points on $U$ by requiring the
intersection with the $\mathcal{Z}$ to be supported on the fibers
above a finite collection $S$ of non-archimedean places of $k$.
(The notion of integral and $S$-integral points on $U$
depends on the choice of integral model $\mathcal{X}$.)

One can speak of the Hasse principle, weak approximation, and strong
approximation for $U$, or more accurately, for
$\mathcal{U}=\mathcal{X}\smallsetminus \mathcal{Z}$.
(For smooth projective varieties, weak and strong approximation are the same.)
Strong approximation concerns the approximation of adelic points by
integral (or $S$-integral) points.
So, the failure of strong approximation can translate as the statement
at a given (non-homogeneous) Diophantine equation has no integral solutions
with prescribed congruence conditions.
We want to illustrate a failure of strong approximation -- concretely
realized as the insolubility of Diophantine equations that admit rational
solutions as well as $v$-adic integer solutions for all $v$ -- which is
explained by the Brauer--Manin obstruction:
$$\mathcal{U}(\mathfrak{o}_k)\subset
\bigl(\prod_{v\nmid\infty} \mathcal{U}(\mathfrak{o}_v)\times
\prod_{v\mid\infty} U(k_v)\bigr)^{\Br(U)},$$
where the set on the right is defined to be the tuples of adelic points,
integral at all non-archimedean places,
whose sum of local invariants is zero with respect
to every element of $\Br(U)$.

The Brauer--Manin obstruction has been extensively studied in the
setting of projective varieties;
see, e.g., \cite{ctks}.
But it has only recently begun to be looked at for open varieties;
see \cite{ctx}.
In the known examples,
$X$ is a quadric surface, $D$ is a conic, and
the relevant Brauer group element is algebraic, i.e.,
lies in the kernel of $\Br(U)\to \Br(U\otimes \bar k)$.
In this paper we use
\emph{transcendental} Brauer group elements.
Transcendental Brauer--Manin obstructions have previously
been exhibited on projective varieties \cite{wittenberg}, \cite{harari}.

We take $U$ to be the complement of a geometrically irreducible smooth
divisor $D$ on $X$.
Then, to have $\Br(X)\ne\Br(U)$ we must have $\dim(X)\ge 2$, so the simplest
examples occur with $X=\PP^2$.
Over an algebraically closed field,
\cite{AM92} relates the Brauer group of the complement of a smooth divisor
to unramified cyclic coverings of the divisor, by means of an
exact sequence:
\begin{equation}
\label{eq:Brseq}
0\to \Br(X\otimes \bar k)\to \Br(U\otimes \bar k)\to
H^1(D\otimes \bar k,\Q/\Z)\to 0.
\end{equation}
We select $D\subset \PP^2$ for which there are
nontrivial unramified coverings defined over $k$,
and for which
(known) constructions of algebras representing ramified Brauer group
elements can be carried out over $k$.
We work over the ground field $k=\Q$.

This note is inspired by lectures of J.-L.~Colliot-Th\'el\`ene on his
joint work with F.~Xu \cite{ctx},
in which they use the Brauer--Manin obstruction
to give a new explanation of the failure of Hasse principle,
exhibited in \cite{BR95} and \cite{SX04}, in the
problem of representing integers by quadratic forms in three variables.

\section{Quartic}

Consider the Diophantine equation
\begin{equation}
\label{eq:qu}
-2x^4-y^4+18z^4=1.
\end{equation}
We claim that
\begin{itemize}
\item[(i)] There are solutions in $p$-adic integers for all $p$.
\item[(ii)] There are solutions in $\Q$.
\item[(iii)] There are no solutions in $\Z$.
\end{itemize}

Indeed, there is the rational solution
$(1/2,0,1/2)$, which is a $p$-adic integer solution for $p\ne 2$.
A $2$-adic integer solution is $(0,\sqrt[4]{17},1)$.
Claim (iii) is more subtle, and uses the Brauer group, as explained
below. 

Our strategy is to relate solutions to \eqref{eq:qu} to points in
the the complement in the projective plane
$$U=\PP^2\smallsetminus D$$
of the divisor $D$ defined by $f(x,y,z)=0$, where
$$f(x,y,z)=-2x^4-y^4+18z^4$$
Then we exhibit a \emph{ramified} Brauer group element $A\in \Br(U)$
corresponding to a degree $2$ cover of $D$ via \eqref{eq:Brseq}.
This forces a congruence condition on integral points on $U$,
incompatible with \eqref{eq:qu}.

The affine variety $U$ has nontrivial geometric fundamental group.
Geometrically,
the universal cover of $U$ is the affine open subset of the $K3$ surface
$f(x,y,z)=t^4$, defined by the
nonvanishing of $t$.
There is an intermediate cover, the degree 2 del Pezzo surface $X$ given by
$$f(x,y,z)=w^2$$
in weighted projective space $\PP(1,1,1,2)$.
The pull-back of $A$ to the del Pezzo surface is \emph{unramified},
i.e., it is the restriction of an element of $\Br(X)$.
The analysis could, then, alternatively be carried out on $X$,
where one would find a Brauer--Manin obstruction to weak approximation,
as in \cite{kt},
that gives us the above-mentioned congruence condition.

In the next section we will exhibit another example
where the obstructing Brauer group element
remains ramified on the universal cover.

By torsor theory, there are finitely many arithmetic twists
$\widetilde{U}_i$ of the geometric universal cover of $U$ such that
$U(\Z)$ is the disjoint union of the images of
the $\widetilde{U}_i(\Z)$.
Obviously if, for each $i$, there is a local
obstruction $\widetilde{U}_i(\Z_{p_i})=\emptyset$,
then $U(\Z)=\emptyset$.
But \eqref{eq:qu} is the defining equation for one of the
$\widetilde{U}_i$, and this has $p$-adic integer points for all $p$,
by claim (i).

By writing the equation for $D$ as
$$(4x^2-y^2)^2 + 2(x^2+2y^2+9z^2)(x^2+2y^2-9z^2)=0$$
we directly see an unramified double cover $\widetilde{D}\to D$
gotten by adjoining the square root of
$(x^2+2y^2+9z^2)/\ell(x,y,z)^2$ to the function field of $D$,
where $\ell$ is an arbitrary linear form.
The recipe of \cite{KRTY} can be carried out
over $\Q$ to produce the following $2$-torsion transcendental
Brauer group element, a quaternion algebra:
\begin{equation}
\label{eq:quaalg1}
(fh,-gh)
\end{equation}
where
\begin{align*}
g(x,y,z)&=-28x^2 - 36xy + 7y^2 + 72z^2,\\
h(x,y,z)&=-25x^2 + 16xy - 22y^2 + 81z^2,
\end{align*}
and where we take the liberty of writing homogeneous functions of
even degree rather than rational functions in
\eqref{eq:quaalg1}.
Essentially, the procedure involves making a coordinate change to
eliminate the $y^4$ term from the quartic equation;
then there is an explicit formula for the quaternion algebra
and the ramification pattern is directly verified to be as desired.
In this case, what needs to be verified is that
$fg$ is a square modulo $h$ and $fh$ is a square modulo $g$;
these verifications can be carried out directly.

The algebra \eqref{eq:quaalg1} defines an element
$A\in \Br(U)$.
This extends to an element of $\Br(\mathcal{U}\otimes \Z[1/2])$,
where $\mathcal{U}=\PP^2_\Z\smallsetminus\mathcal{D}$.
Hence, $A$ is unramified at all points of
$\mathcal{U}(\Z_p)$ for $p\ge 3$.
According to a $2$-adic analysis,
at $\Z_2$-points $(X:Y:Z)$ of $\PP^2$ with
\begin{equation}
\label{eq:mod2cond}
X\equiv 0 \pmod 2,\qquad
Y\equiv 1\pmod 2,\qquad
Z\equiv 1\pmod 2,
\end{equation}
we have
\begin{align*}
f(X,Y,Z)h(X,Y,Z)&\equiv 3\pmod 4,\\
-g(X,Y,Z)h(X,Y,Z)&\equiv 3\pmod 4,
\end{align*}
hence $A$ is \emph{ramified} at all such
$2$-adic points.
We leave it to the reader to check that
$A$ is unramified at real points of $U$.
So, for any integers $X$, $Y$, $Z$ satisfying \eqref{eq:mod2cond}
there must exist a prime divisor $p\ge 3$ of $f(X,Y,Z)$.
Examining reductions modulo $16$,
any integer solution to
\eqref{eq:qu} would have to satisfy \eqref{eq:mod2cond},
hence \eqref{eq:qu} has no integer solutions.
Notice that there is an obvious integer point $(0:1:0)$ on $\mathcal{U}$.
So, this example constitutes
a Brauer--Manin obstruction to strong approximation on $\mathcal{U}$.

\section{Cubic}

Consider the Diophantine equation
\begin{equation}
\label{eq:cu}
y^2z-(4x-z)(16x^2+20xz+7z^2)=1.
\end{equation}
We claim, again, that there are solutions in $\Z_p$ for all $p$,
in $\Q$, but not in $\Z$.
A rational solution is $(1/4,1,1)$, and a $2$-adic integer solution
is $(0,0,\sqrt[3]{1/7})$.

Define $D$ to be the divisor $f(x,y,z)=0$, where
$$f(x,y,z)=y^2z-(4x-z)(16x^2+20xz+7z^2).$$
So $D$ is an elliptic curve, and one easily computes its group structure
$D(\Q)\cong \Z/2\Z$.
Set $U=\PP^2\smallsetminus D$.
There is an element $A\in \Br(U)$
corresponding to the (unique) unramified degree $2$ cover of $D$ defined
over $\Q$, by \cite{Jacob}:
the quaternion algebra
\begin{equation}
\label{eqn:anotherquaalg}
(y^2z^2-(4x-z)(16x^2+20xz+7z^2)z, (4x-z)z).
\end{equation}
By an analysis similar to that above, $A$ is
ramified at $2$-adic points satisfying
\begin{equation}
\label{eqn:second2adiccond}
Y\equiv 0\pmod 2\qquad\text{and}\qquad Z\equiv 1\pmod 2,
\end{equation}
and is unramified at real points and $\Z_p$-points of $U$,
for $p\ge 3$. A quick analysis modulo 2 reveals that
any integer solution $(X:Y:Z)$ to \eqref{eq:cu} must satisfy
\eqref{eqn:second2adiccond},
so as above we conclude that \eqref{eq:cu} has no
integer solutions.
In fact, $\mathcal{U}(\Z)=\emptyset$, where
$\mathcal{U}=\PP^2_\Z\smallsetminus\mathcal{D}$,
i.e., in this example we have a
Brauer--Manin obstruction to the Hasse principle over $\Z$.
Also, since $U$ has geometric fundamental group $\Z/3\Z$,
the class of \eqref{eqn:anotherquaalg}
in the Brauer group remains transcendental on the universal cover.

\end{document}